\documentclass[11pt]{article}
\usepackage{amssymb}
\def\qed{\hbox
to 0pt{}\hfill$\rlap{$\sqcap$}\sqcup$}

\newcommand{\dem}{\noindent{\bf Proof. }}
\newcommand{\fg}{{\mathfrak g}}
\newcommand{\fa}{{\mathfrak a}}
\newcommand{\fh}{{\mathfrak h}}
\newcommand{\ad}{\mbox{\rm ad}}
\newcommand{\ada}{\mbox{\rm ad}_\fa}
\newcommand{\adg}{\mbox{\rm ad}_\fg}

\newtheorem{theor}{Theorem}[section]

\newtheorem{prop}[theor]{Proposition}
\newtheorem{coro}[theor]{Corollary}
\newtheorem{lema}[theor]{Lemma}
\newtheorem{rema}{Remark}

\newtheorem{defi}{Definition}[section]

\let \a = \alpha
\let \d = \delta

\let\o = \omega

\let \g = \gamma 
\let \t = \theta

\let \l = \lambda

\let \f = \mathfrak
\let \x = \backslash

\let \n = \noindent

\author{Ignacio Bajo, Sa\"{\i}d Benayadi and Alberto Medina}
\title{Symplectic structures on quadratic Lie algebras}
\date{ }
\begin{document}
\maketitle

\begin{abstract}
We study quadratic Lie algebras over a field $\Bbb K$ of null characteristic which admit, at the
same time, a symplectic structure. We see that if  $\Bbb K$ is algebraically closed every such
Lie algebra may be constructed as the T$^*$-extension of a nilpotent algebra admitting an invertible
derivation and also as the double extension of another quadratic symplectic Lie algebra by the
one-dimensional Lie algebra.  Finally,   we  prove
that every symplectic quadratic Lie algebra is a special symplectic Manin algebra and we give an inductive
classification in terms of symplectic quadratic double
extensions.  
\end{abstract}

\section*{Introduction}

Lie groups which admit a bi-invariant pseudo-Riemannian metric and also a left-invariant symplectic form are nilpotent Lie groups and their geometry (and, consequently, that of their associated homogeneous spaces) is very rich. In particular, they carry two left-invariant affine structures: one  defined by the symplectic   form (which is well-known) and another which is compatible with a left-invariant pseudo-Riemannian metric. Moreover, if the symplectic form is viewed as a solution $r$ of the classical Yang-Baxter equation, then the Poisson-Lie tensor $\pi =r^+-r^-$ and the geometry of the double Lie groups $D(r)$ can be nicely described \cite{diatta}. As we will see below, a great number of such Lie groups may be constructed and, hence, a large class of symplectic symmetric  spaces is found. 

The Lie algebra of one of those Lie groups turns to be a {\it quadratic symplectic Lie algebra}, this is to say, a Lie algebra endowed with both an invariant
non-degenerate symmetric bilinear form and a non-degenerate
2-cocycle of its scalar cohomology. In this paper we study the structure of these nilpotent Lie algebras over a field $\Bbb K$ of characteristic zero and give some results which provide a method for their inductive classification. 
The main tools used for our purposes are T$^*$-extensions and quadratic double extensions. 

The notion of T$^*$-extension of a Lie algebra was introduced by M. Bordemann in \cite{borde}, where
it is proved that every solvable quadratic Lie algebra over an algebraically closed field of
characteristic 0 is either a T$^*$-extension or a non-degenerate ideal of codimension 1 in a
T$^*$-extension of some Lie algebra. In general, a T$^*$-extension of a  Lie algebra need not admit
a non-degenerate scalar 2-cocycle, even if the extended algebra is nilpotent. Thus, additional
properties should be imposed to a Lie algebra $\fg$ to ensure that a T$^*$-extension
$T^*_{\theta}\fg$ might be furnished with a symplectic structure. This is done in section 2, where
we show that if $\Bbb K$ is algebraically closed, then every  quadratic symplectic Lie algebra  is a T$^*$-extension of a Lie algebra
which admit an invertible derivation and we give necessary and sufficient conditions on $\fg$ and
on the cocycle $\theta$ used in the construction of the T$^*$-extension $T^*_{\theta}\fg$ to assure
that the extended algebra admits an skew-symmetrical derivation and, hence, a symplectic structure.
We use these results to give a complete classification of complex quadratic Lie algebras of dimension less than or equal to 8 which admit a symplectic structure.

Every $n$-dimensional solvable quadratic Lie algebra may be obtained from a quadratic Lie algebra
of dimension $n-2$ by a central extension by a one-dimensional algebra and a semi-direct product
by another one-dimensional algebra. This method of construction, known as quadratic double extension, was
introduced and developped for the first time in \cite{m-r-travaux}, \cite{m-r85}. Later, Dardi\'e and Medina addapted and generalised this method to study certain symplectic Lie groups \cite{dard}. In A. Aubert's unpublished thesis 
\cite{aubert}, these methods of (generalised) symplectic double extension are used to study quadratic symplectic Lie algebras. A restatement of some of her results using quadratic   double extension instead of symplectic double extension is given in section 3 since it will be useful in the last section of the paper, where we study Manin algebras admitting an adapted symplectic form, which we will call {\it special symplectic Manin algebras}. This study is motivated by the fact that every symplectic quadratic Lie algebra over an algebraically closed field is actually a special symplectic
Manin algebra. We give a description of these Lie algebras in term of quadratic double extension and obtain that, if the field is algebraically closed, every such special  symplectic Manin algebra   may be
obtained from a two-dimensional symplectic Manin algebra by a sequence of  quadratic 
double extensions by the one-dimensional Lie algebra where the algebra obtained in each step is also a special symplectic Manin algebra.

\section{Definitions and preliminary results}

All through the paper $\Bbb K$ will denote a commutative field of characteristic 0.
We begin with the following:

\begin{defi}{\rm Let $\fg$ be a Lie algebra over $\Bbb K$.
\begin{enumerate}
\item[i)]
 We say that
$(\fg,B)$ is a {\it quadratic Lie algebra}  if $B$ is a non-degenerate  symmetric bilinear form on
$\fg$ such that
$B([x,y],z)=B(x,[y,z])$ for all $x,y,z\in\fg$. In that case, we will say that $B$
is an {\it invariant scalar product} on $\fg$.

A quadratic Lie algebra $(\fg,B)$ is said to be {\it reducible} (or $B$-reducible) if it admits an ideal $\mathfrak I$ such that the restriction of $B$ to ${\mathfrak I}\times{\mathfrak I}$ is non-degenerate. We will say that $(\fg,B)$ is  {\it irreducible} otherwise.
\item[ii)]
We say that
$(\fg,\o)$ is a {\it symplectic Lie algebra}  if $\o$ is a non-degenerate skew-symmetric bilinear
form on $\fg$ and
$\omega ([x,y],z)+\omega ([y,z],x)+\omega ([z,x],y)=0$ holds for all $x,y,z\in\fg$, this is to
say, $\o$ is a non-degenerate 2-cocycle for the scalar cohomology of $\fg$. 
Note that in such case, $\fg$  must be even-dimensional. We will then call  $\o$ a {\it symplectic
structure} on $\fg$. 
\item[iii)] We will say that $(\fg,B,\o)$ is a {\it quadratic symplectic Lie algebra} if
$(\fg,B)$ is  quadratic and $(\fg,\o)$ is  symplectic.
\end{enumerate}}\end{defi}
\vspace*{0.15cm}

\begin{lema}\label{existder}
Let $({\mathfrak g},B)$ be a quadratic Lie algebra over ${\Bbb K}$. A symplectic structure $\omega$ may be defined on $\mathfrak g$ if and only if
there exists a skew-symmetric invertible derivation $D$ of $({\mathfrak g},B)$.
\end{lema}
\dem It is a straightforward calculation considering
$\omega (x,y)=B(Dx,y)$
for all $x,y\in\fg$. \qed

\vspace*{0.25cm}

{From} now on, given a quadratic Lie algebra $(\fg ,B)$,  we will denote by $\mbox{Der}_a(\fg ,B)$
the Lie algebra of skew-symmetric derivations of $(\fg ,B)$.

\begin{rema}\em
\begin{enumerate}\item[$\,$]
\item Since the Lie algebra $\mbox{Der}_a(\fg ,B)$ is algebraic, the existence of an invertible skew-symmetric derivation of $({\mathfrak g},B)$ is equivalent to the existence of an invertible  semisimple skew-symmetric derivation of $({\mathfrak g},B)$.
\item Clearly, every quadratic Lie algebra admitting a symplectic structure (and hence an invertible derivation) must be nilpotent
\cite{jacob}.
\item It should be noticed that under the assumptions of the lemma above, 
the skew-symmetric derivation $D$ of $(\fg ,B)$ is also skew-symmetric with respect to the
symplectic form $\omega$ since for all $x,y\in\fg$ we get
$$\omega (Dx,y)=B(D^2x,y)=-B(Dx,Dy)=-\omega (x,Dy).$$

One might thus think that every symplectic Lie algebra $(\fg ,\omega)$ admitting an  invertible derivation which is skew-symmetric for $\omega$ carries a quadratic structure; but this is not the case. For example,  the four dimensional Lie algebra $\fg ={\Bbb R}\mbox{-span}\{x_1,x_2,x_3,x_4\} $ defined by the only non-trivial bracket $[x_1,x_2]=x_3$ does not admit any quadratic structure since $[\fg ,\fg]$ is one dimensional while the centre of $\fg$ has dimension two. However, the skew-symmetric bilinear form on $\fg$ given by $\omega (x_1,x_4)=\omega (x_2,x_3)=1$ provides a symplectic structure on $\fg$ and the linear endomorphism of $\fg$ given
by
$D(x_1)=2x_1 ,\, D(x_2)=-x_2 ,\, D(x_3)=x_3 ,\, D(x_4)=-2x_4,$
is a skew-symmetric derivation of $(\fg ,\omega)$.
\end{enumerate}
\end{rema}

\bigskip

The following example shows that, starting from an arbitrary Lie algebra, one can construct infinitely many quadratic symplectic Lie algebras.

\medskip

\noindent {\bf Example 1} $\,$
Let  $\f g$ be a Lie algebra and $n\in {\Bbb N},\ n>1.$ If we consider the
non-unitary  associative  algebra ${\cal A}_n= {t{\Bbb K}[t]/t^n{\Bbb K}[t]},$   Then the vector space 
${{\cal L}_n}= {\f g}\otimes {\cal A}_n$ with the bracket
\begin{eqnarray*}
&[x\otimes {\bar{t^p}},y\otimes {\bar{t^q}}]= [x,y]_{\f g}\otimes {\bar{t^{p+q}}}, \quad x, y \in
{\f g}, \,\, p,q \in {\Bbb N}\setminus \{0\},& 
\end{eqnarray*}
is a nilpotent Lie algebra.
The endomorphism $D$   of ${\cal L}_n$ defined by  
$D(x\otimes {\bar{t^p}})= p(x\otimes {\bar{t^p}}),$ for all $x \in
{\f g}, \ p \in \{1,\dots,n-1\},$
is an invertible derivation of ${\cal L}_n.$

Now, the vector space ${\tilde {\cal L}_n}= {\cal L}_n\oplus({\cal L}_n )^*$ with the bracket defined by 
$$[X+f,Y+g]= [X,Y]_{{\cal L}_n} - g\circ\ad _{{\cal L}_n}(X)  +f\circ\ad _{{\cal L}_n}(Y),$$ for
$X,Y \in {\cal L}_n, \ f, g \in ({\cal L}_n )^*,$ 
and the bilinear form   
$ B(X+f,Y+g)= f(Y) + g(X),$
is a quadratic Lie algebra. Further an
invertible skew-symmetric derivation of  ${\tilde {\cal L}_n}$ may be defined by 
$ {\tilde D}(X+f)= D(X) + D^*(f),$  $\forall X \in {\f g},
~ f \in {\f g}^*,$ where $D^*(f)= -f\circ D, $  and hence the quadratic algebra $({\tilde {\cal
L}_n},B)$ admits a symplectic structure.

\vspace*{0.5cm}

It is well-known that if $(\fg ,\omega)$ is a symplectic Lie algebra then the product $x.y$ on $\fg$ defined by 
$\omega (x.y,z)=-\omega (y,[x,z])$
induces on every Lie group $G$ with Lie algebra $\fg$ a flat and torsion-free invariant connection by the formula
$\nabla ^\omega_xy=x.y,$ for $x,y\in\fg.$
Further, in \cite{m-r91} it is proved that if $\fg$ also admits a quadratic structure then one may
define a left-invariant pseudo-Riemannian metric $\langle .,.\rangle$ on $G$ whose Levi-Civita
connection is precisely the flat connection $\nabla ^\omega$. A slight generalisation of such situation may be done since the result remains valid for non-symplectic Lie algebras which admit a (not necessarily skew-symmetric) invertible derivation. In order to prove this, we recall the following well-known result:

\begin{lema} If $\fg$ is a Lie algebra and $D$ is an invertible derivation of $\fg$ then the product
$x.y=D^{-1}[x,Dy]$ for $x,y\in\fg,$
defines a Lie-admissible complete left-symmetric structure on $\fg.$ 
Further, $D$ is also a derivation of the left-symmetric product.
\end{lema}

\begin{prop}
Let $\fg$ be a Lie algebra, $D$ an invertible derivation of $\fg$ and let $(\fg , .)$ be the left-symmetric structure defined by $D$. The algebra
$\fg$ admits a quadratic structure
if and only if there exists a non-degenerate symmetric bilinear form 
$\langle .,.\rangle$ on $\fg$ such that its Levi-Civita connection is given by
$\nabla _xy=x.y,$ for all $ x,y\in\fg.$
\end{prop}
\dem The proof is almost immediate if we define $\langle .,.\rangle$ or, conversely,  the invariant scalar product $B$ according to
 $\langle x,y\rangle =B(Dx,Dy),$
for all
$x,y\in\fg.$\qed

\vspace{0.25cm}

The following corollary is the obvious geometric analog of the proposition above.

\begin{coro}
Let $G$ be a Lie group and suppose that its Lie algebra $\mathfrak g$ admits an invertible derivation $D$. The following are
equivalent:
\begin{enumerate}
\item[(i)] $G$ admits a bi-invariant pseudo-Riemannian metric.
\item[(ii)] $G$ admits a left-invariant flat pseudo-Riemannian metric whose Levi-Civita connection
is given by left-translation of the left-symmetric product induced by $D$ on
$\fg$.
\end{enumerate}\end{coro}

\smallskip

\begin{rema} {\em Quadratic Lie algebras admitting an invertible skew-symmetric derivation provide nice 
examples of solutions of both classical and modified Yang-Baxter equations. We recall that for a quadratic Lie algebra $(\fg,B)$ a linear endomorphism $R$ of $\fg$ which is skew-symmetric with respecto to $B$ is said to
satisfy the classical Yang-Baxter equation  if 
$[Rx,Ry]-R[Rx,y]-R[x,Ry]=0$
holds for every $x,y\in\fg$.
If, instead, for every $x,y\in\fg$ the mapping $R$ verifies
$[Rx,Ry]-R[Rx,y]-R[x,Ry]+[x,y]=0$
one says that $R$ is a solution of the modified Yang-Baxter equation. A straightforward calculation shows that if $D$ is an invertible skew-symmetric derivation of $(\fg,B)$ then $R=D^{-1}$ is a solution of the classical Yang-Baxter equation. Further, one easily sees that in the complex case a quadratic Lie algebra with an invertible skew-symmetric derivation defines naturally a Manin triple and, as a consequence, a solution of the modified Yang-Baxter equation \cite{al-mi}, \cite{CP}, \cite{m-r88}.}
\end{rema}

\noindent {\bf Notation } In all the sections below, the symbol $\oplus$ will be frequently used. In principle, unless other thing is stated, it will only denote direct sum of vector spaces.

\section{Quadratic Lie algebras with symplectic structures and the notion of $T^*$-extension}

In \cite{borde}, M. Bordemann defines the following notion:

\begin{defi}\label{borde} {\rm \cite{borde}\, Let $\fa$ be a Lie algebra over a commutative 
field and consider its dual space $\fa^*$. Consider a 2-cocycle $\theta
:\fa\times\fa\rightarrow\fa^*$ and define on the vector space $T^*_\theta\fa=\fa\oplus\fa^*$ the
following bracket:
$$[x+f,y+g]=[x,y]_\fa +\theta (x,y) +f\circ\ada (y)-g\circ\ada (x)\, ,\quad x,y\in\fa,\, 
f,g\in\fa^*.$$
The pair $(T^*_\theta\fa ,[.,.])$ is then a Lie algebra called the {\it $T^*$-extension of
the Lie algebra $\fa$ by means of $\theta$}.

 If, further, the 2-cocycle $\theta$ verifies the {\it cyclic} condition
$\theta (x,y)(z)=\theta (y,z)(x)$ for all $
x,y,z\in\fa,$ then the symmetric bilinear form $B$ on $T^*_\theta\fa$
given by
$B(x+f,y+g)=f(x)+g(y)$ for $x,y\in\fa,\, f,g\in\fa^*,$
defines a quadratic structure on $T^*_\theta\fa$.
}\end{defi}

Note that the condition of cyclicity for $\theta$ proves that $\theta$ defines a 3-cocycle of the scalar cohomology. Actually, the natural mapping between the set of cyclic 2-cocyles and $Z^3(\fa ,\Bbb K)$ given by $\tilde{\theta}(x,y,z)= 
\theta (x,y)(z)$ for all $x,y,z$ is an isomorphism. 

\medskip

The case $\theta =0$ corresponds to the double extension of the null algebra by the algebra $\fa$ \cite{m-r85}. This algebra, which we will call the {\it trivial $T^*$-extension of $\fa$}, is isometrically isomorphic to  $T^*_\theta\fa$ if $\tilde{\theta}$ is a coboundary \cite{borde}.

\medskip 

Since we will mainly deal with even-dimensional Lie algebras the following particular case
of
 \cite[Cor. 3.3]{borde} will be sufficient for our purposes:

\begin{prop}\label{martin} Suppose that $\Bbb K$ is algebraically closed and  let $(\fg,B')$ be a  quadratic even-dimensional Lie algebra over $\Bbb K$. 
If $\fg$ is solvable then $(\fg,B')$ is isometrically isomorphic to a quadratic $T^*$-extension
$(T^*_\theta\fa ,B)$, where $\fa$ is isomorphic to the quotient algebra of $\fg$ by a completely isotropic ideal.
\end{prop}

 From the proposition above, one obviously gets that every  Lie algebra admitting both a 
quadratic and a symplectic structure is isometrically isomorphic to some $T^*$-extension. However
it is clear that not every $T^*$-extension of a nilponent Lie algebra should admit a symplectic
structure. Our main interest is to give necessary and sufficient conditions for a quadratic Lie
algebra endowed with a symplectic structure to be a $T^*$-extension of some subalgebra. This is
done in the two following propositions.

\begin{prop} \label{T*up} Let $\fa$ be a Lie algebra admitting an invertible derivation $D$. 
Consider a  cyclic 2-cocycle $\theta\in Z^2({\mathfrak a},{\mathfrak a^*})$   and
define
$$\Theta (x,y,z)=\theta (Dx,y)z+\theta (Dy,z)x+\theta (Dz,x)y\, ,\quad x,y,z\in {\mathfrak a}.$$

If $\Theta$ is a 3-coboundary for the scalar cohomology of $\fa$ then the
quadratic Lie algebra $T^*_\theta \mathfrak a$ admits a  symplectic structure.
\end{prop}

 \dem Let $B$ be the quadratic form on  $T^*_\theta \mathfrak a$ defined in Definition \ref{borde}. 
By Lemma \ref{existder} it suffices to prove the existence of an invertible skew-symmetric
derivation of $(T^*_\theta \mathfrak a,B).$

Let us consider $F:\fa\wedge\fa\rightarrow\Bbb K$ be such that $\Theta =dF$ and let $H:\mathfrak a\rightarrow\mathfrak a^*$ be the mapping defined
by
$B(Hx,y)=F(x,y)$ for all $x,y\in\mathfrak a.$
Now, define $\overline{D}: T^*_\omega {\mathfrak a}\rightarrow T^*_\omega {\mathfrak a}$  by
$\overline{D}(x+f)=Dx -Hx-f\circ D$ for all $x\in {\mathfrak a},
f\in {\mathfrak a}^*.$
It is straightforward to see that $\overline{D}$ is invertible, since $D$ is so, and a direct calculation shows that
$\overline{D}$ is also a skew-symmetric with respect to $B.$

Further, since $D$ is a derivation of $\fa$, we
get
\begin{eqnarray}\nonumber
& & [\overline{D}(x+f),y+g]+[x+f,\overline{D}(y+g)]-
\overline{D}[x+f,y+g]=\\\label{prueder} & &
\quad \theta (x,y)\circ D+\theta (Dx,y)+\theta (x,Dy)+ H[x,y]_\fa-Hx\circ \ada (y)+Hy\circ\ada
(x).\end{eqnarray} But since $\Theta =dF$ and $F(x,y)=B(Hx,y)=Hx(y)$ for all
$x,y\in\fa$, we have
\begin{eqnarray*}&& \theta (x,y)Dz+\theta (Dx,y)z+\theta (x,Dy)z+ H[x,y]_\fa (z)-Hx([y,z]_\fa )+Hy([x,z]_\fa
)=\\ & & \quad\Theta(x,y,z)+F([x,y]_\fa,z)-F(x,[y,z]_\fa )+F(y,[x,z_\fa ])=\\&
&\quad\Theta(x,y,z)-dF(x,y,z)=0
\end{eqnarray*}
and, hence, (\ref{prueder}) vanishes.\qed

\bigskip

If the field $\Bbb K$ is algebraically closed, one also has the reciprocal:

\begin{theor}\label{T*down} Let $({\fg} ,B)$ be a quadratic  Lie algebra over an algebraically closed field $\Bbb K$ which admits a skew-symmetric  invertible derivation $\overline{D}$.

There exist a Lie algebra $\fa$, an invertible derivation $D$ of $\fa$ and a cyclic $\theta\in Z^2({\fa},{\fa ^*})$ such
that 
 ${\mathfrak g}=T^*_\theta{\mathfrak a},$
and the map $\Theta$ defined by $\Theta (x,y,z)=\theta (Dx,y)z+\theta (Dy,z)x+\theta (Dz,x)y$, for all$x,y,z\in\fa$ is a 3-coboundary for the scalar cohomology of
$\fa$. 
\end{theor}
\dem Let us consider the semidirect sum of  Lie algebras 
${\mathfrak L}=\ad (\fg)\oplus{\Bbb K}\overline{D}$. Since $\fg$ is nilpotent, the Lie algebra $\mathfrak
L$ is obviously solvable. Thus, according to Lemma 3.2 in \cite{borde}, we may find a maximal
isotropic (with respect to the quadratic form $B$) ideal $\mathfrak I$ of $\fg$ which is also stable by
the derivation $\overline{D}.$ Now, if $\fa =\fg /{\mathfrak I}$ then $\fa ^*$ is isomorphic to $\mathfrak
I$ and there exists a cyclic $\theta\in Z^2({\fa},{\fa ^*})$ such that $(\fg ,B)$ is isometrically
isomorphic to $T^*_\theta\fa$ \cite[Corollary 3.1]{borde}. Furhter, $\fa ^*=\mathfrak I$ is stable by
$\overline{D}$  and hence, there exist linear mappings
$D_{11}:\fa\rightarrow\fa$, $D_{21}:\fa\rightarrow\fa ^*$ and $D_{22}:\fa ^*\rightarrow\fa ^*$ such
that
$\overline{D}(x+f)=D_{11}x+D_{21}x+D_{22}f$ holds
for every $x\in\fa$, $f\in\fa ^*$. Clearly, $D_{11}$ and $D_{22}$ must be invertible since $\overline{D}$ is so. The skew-symmetry of $\overline{D}$ is equivalent to the 
conditions $D_{22}f=-f\circ D_{11}$, for all $f\in\fa^*$ and $B(D_{21}x,y)=-B(D_{21}y,x)$ for $x,y\in\fa$
Let $H=-D_{21}$ and $D=D_{11}$. Since $\overline{D}$ is a derivation, we get:
\begin{eqnarray*}
& & 0=[\overline{D}x,y]+[x,\overline{D}y]-
\overline{D}[x,y]=[Dx,y]_\fa+[x,Dy]_\fa-D[x,y]_\fa+\\ & & \quad\quad
H[x,y]_\fa-Hx\circ\ada (y)+Hy\ada (x)+\theta (x,y)\circ D +\theta (Dx,y)+\theta (x,Dy),
\end{eqnarray*}
for $x,y\in\fa$, which shows that $D$ is a derivation of $\fa$ and that if 
$$\Theta (x,y,z)=\theta (Dx,y)z+\theta (Dy,z)x+\theta (Dz,x)y$$ and
$F$ is the skew-symmetric bilinear form on $\fa$ defined by $F(x,y)=-B(Hx,y)=-Hx (y),$ then we
have
$$\Theta (x,y,z)+F([x,y]_\fa ,z)-F(x,[y,z]_\fa)-F(y,[z,x]_\fa)=\Theta
(x,y,z)-dF(x,y,z)=0$$ for all $x,y,z\in\fa$, which finishes the proof.\qed

\medskip

\begin{rema}{\em \begin{enumerate}
\item[$\,$]
\item It is interesting to point out that if the derivation $\overline{D}$ is semisimple then the derivation $D$ is also semisimple. Actually, if we choose a basis of the completely isotropic ideal ${\mathfrak I}$ composed of eigenvectors of $\overline{D}$ then each element in its dual basis (with respect to $B$) is an eigenvector of $D$.
\item
Theorem \ref{T*down} does not hold in general in the  case of a non-algebraically closed field. For example, an even-dimensional abelian Lie algebra $\mathfrak g$ over $\Bbb R$ with a definite positive bilinear form is obviously a quadratic Lie algebra which admits an invertible skew-symmetric derivation. However, it cannot be a  $T^*$-extension since there are no isotropic subspaces.
\end{enumerate}}
\end{rema}

\medskip

In \cite{jacob} it is proved that if a Lie algebra admits an invertible derivation, then it must be nilpotent. There are however many nilpotent Lie algebras whose derivations are all singular. The following result gives a characterization of Lie algebras admitting such a derivation. Note that the result is valid for an arbitrary base field of characteristic zero (not necessarily algebraically closed). We recall that in a symplectic Lie algebra $(\fg,\omega)$ an ideal is called {\it lagrangian} if and only if it coincides with its orthogonal with respect to the form $\omega$. 

\begin{prop}
Let $\Bbb K$ be a field of characteristic zero and let $\fa$ be a Lie algebra over $\Bbb K$. 

There exists an invertible derivation of $\fa$ if and only if $\fa$ is isomorphic to the quotient Lie algebra 
$\fg /{\mathfrak I}$ of a quadratic symplectic Lie algebra $(\fg, B, \omega)$ by a lagrangian and completely isotropic ideal $\mathfrak I$.
\end{prop}
\dem If $\fa$ admits an invertible derivation then the Lie algebra $\fg=T^*_0\fa$ obtained by $T^*$-extension by the null cocycle $\theta =0$ is, according to Proposition \ref{T*up}, a quadratic symplectic Lie algebra and  ${\mathfrak I}=\fa^*$ is a lagrangian, completely isotropic ideal of $\fg$.

Conversely, suppose that $\fa$ is isomorphic to 
$\fg /{\mathfrak I}$ where $(\fg, B, \omega)$ is a quadratic symplectic Lie algebra and $\mathfrak I $ is a lagrangian completely isotropic ideal of $\fg$. According to \cite[Corollary 3.1]{borde}, $\fg$ is isometrically isomorphic to $T^*_\theta (\fg / {\mathfrak I})=T^*_\theta\fa$ since $\mathfrak I$ is completely isotropic. Let $\overline{D}\in \mbox{\rm Der}_a( \fg,B)$ be the invertible derivation such that 
$\omega (X,Y)=B(\overline{D}X,Y)$ for all $X,Y\in\fg$. Clearly, $\omega ({\mathfrak I},{\mathfrak I})=\{0\}$ implies
that $\overline{D}({\mathfrak I})\subset {\mathfrak I}^\perp={\mathfrak I}$. Now, since, ${\mathfrak I}$  stable by $\overline{D}$, the same arguments used in the proof of Theorem \ref{T*down} prove that the projection of $\overline{D}_{\vert\fa}$ to $\fa$ provides a non-singular derivation of $\fa$.

\section{Classification of  quadratic  symplectic complex Lie algebras up to dimension 8}

The results of the preceding section allow the classification of Lie algebras which admit a quadratic and a symplectic structure of a given dimension by the calculation of the nonisomorphic $T^*$-extensions of all the Lie algebras whose dimension is exactly the half. We  first reduce the classification problem to the study of $T^*$-extensions of Lie algebras which do not admit an abelian direct summand. The main result is the following:

\begin{prop}\label{red}
Let $\fa=\fh\oplus {\Bbb C}$ be a direct sum of Lie algebras and let $\theta\in Z^2(\fa,\fa^*)$ be a cyclic cocycle. 
\begin{enumerate}
\item[(i)] If the $T^*$-extension $\fg=T^*_\theta (\fa)$ is irreducible, then there exists a Lie algebra $\fa _1$ and a cyclic $\theta _1\in Z^2(\fa,\fa^*)$ such that
\begin{eqnarray*}
 \fg=T^*_{\theta _1} (\fa_1)\, ,&  \dim {\mathfrak z}(\fa _1)\le \dim {\mathfrak z}(\fa)\, , &
\dim \Bigl([\fa _1,\fa_1]_{\fa_1}\cap {\mathfrak z}(\fa _1)\Bigl)=\dim \Bigl([\fa,\fa]_{\fa}\cap {\mathfrak z}(\fa)\Bigl)+1.\end{eqnarray*}
\item[(ii)] If, further, $\overline{D}$ is a semisimple skew-symmetric invertible derivation of $\fg$ leaving $\fa^*$ invariant then $\mathfrak a_1$ may be chosen such that $\mathfrak a_1^*$ is also stable by $\overline{D}.$ 
\end{enumerate}
\end{prop}
\dem In order to prove (i), let ${\cal B}=\{x_1,x_2,\cdots ,x_n,e\}$ be a basis of $\fa$ where $x_i\in\fh$, for all $i\le n$, and $e\in\mathfrak{z}(\fa)$ such that $e\not\in\mathfrak{z}(\fh)$ and consider its dual basis ${\cal B}^*=\{x_1^*,x_2^*,\cdots ,x_n^*,e^*\}$. There must exist two elements $x_k,x_l\in\fa$ such that $\theta (x_k,x_l)(e)\neq \{0\}$ since, otherwise, $\fg=T^*_\theta (\fa)$ would split as a orthogonal sum $T^*_\theta (\fa)=T^*_{\overline{\theta}} (\fh)\oplus\mbox
{${\Bbb C}$-span}\{e,e^*\},$ where $\overline{\theta}$ is the restriction of $\theta$ to $\fh\times\fh$.  Recall that the center of $\fg$ is given by:
$${\mathfrak z}(\fg)=\{x+f\in \fa\oplus\fa^*\, ;\, f([y,z]_\fa)=0, [x,y]_\fa=0 \mbox{ and } \theta(x,y)=0, \mbox{ for all }y,z\in\fa\}.$$
Therefore, $e\not\in{\mathfrak z}(\fg)$ and actually ${\mathfrak z}(\fg)\subset {\mathfrak z}(\fh)\oplus\fa^*$, which implies that $e^*\in {\mathfrak z}(\fg)^\perp=[\fg,\fg]$.

Now, consider de linear subspace of $T^*_\theta (\fa)$ given by
${\mathfrak I}=\mbox{${\Bbb C}$-span}\{x_1^*,x_2^*,\cdots ,x_n^*,e\}.$
It is clear that $\mathfrak I$ is a subalgebra of $\fg$, because $e$ is a central element of $\fa$, and that it is completely isotropic. Notice that $\mathfrak I$ actually contains $\fh^*$ and that $\theta (x,e)(e)=0$ for every $x\in\fa$. Therefore, $\mathfrak I$  is an ideal of $\fg$ since $[e^*,x_i^*]=[e^*,e]=0$, $[x_i,x_j^*]=x_j^*\circ\ad _\fa (x_i)\in \fh^*$ and $[x_i,e] =\theta (x_i,e)\in \fh^*$, for all $i,j\le n$. 
Thus, according to Corollary 3.1 in \cite{borde}, the algebra $\fg$ is isometrically isomorphic to a $T^*$-extension
$T^*_{\theta _1} (\fa_1)$ of $\fa_1=\fg /{\mathfrak I}$. Note that $\fa_1$ is spanned by the set 
$\{p({x}_1),p({x}_2),\cdots ,p({x}_n),p({e}^*)\}$, where $p:\fg\rightarrow\fg /{\mathfrak I}$ denotes the canonical projection. The restriction $p_{\vert \fh\oplus{\Bbb C}e^*}$ is clearly injective and  the brackets in $\fa _1$ are given by:
$$[p({x}_i),p({x}_j)]_{\fa_1}=p([x_i,x_j]_{\fa})+\theta (x_i,x_j)(e)p(e^*)\, ,[p({x}_i),p({e^*})]_{\fa_1}=0\, ,\quad 1\le i,j\le n.$$
Recall that $e^*\in[\fg,\fg]$ but it cannot be of the form $f\circ \ad_\fa(x)$ for $x\in\fa, f\in\fa^*$ because, otherwise, we would have $e^*(e)=f\circ \ad_\fa(x)(e)=0$, an absurd. Thus, we have that $e^*\in[\fh,\fh]$ and therefore $p({e^*})\in[\fa_1,\fa_1].$ This proves that 
$[\fa_1,\fa_1]_{\fa_1}=p\Bigl([\fh,\fh)]_\fa\Bigl)\oplus{\Bbb C}p({e^*})$ and further
$[\fa_1,\fa_1]_{\fa_1}\cap{\mathfrak z}(\fa_1)=p\Bigl([\fh,\fh]_\fa\cap{\mathfrak z}(\fh)\Bigl)\oplus{\Bbb C}p({e^*})$. Since $[\fh,\fh]_\fa\cap{\mathfrak z}(\fh)=[\fa,\fa]_\fa\cap{\mathfrak z}(\fa)$ this proves that $\dim \Bigl([\fa _1,\fa_1]_{\fa_1}\cap {\mathfrak z}(\fa _1)\Bigl)=\dim \Bigl([\fa,\fa]_{\fa}\cap {\mathfrak z}(\fa)\Bigl)+1$. Finally, ${\mathfrak z}(\fa)={\mathfrak z}(\fh)\oplus{\Bbb C}{e}$ and one has
${\mathfrak z}(\fa_1)=\{p(x)\, :\, x\in{\mathfrak z}(\fh) \mbox{ and }\theta (x,e)=0\}\oplus{\Bbb C}p(e^*),$
 which  clearly shows that  $\dim {\mathfrak z}(\fa_1)\le \dim {\mathfrak z}(\fa)$.

For the proof of (ii), if  $\overline{D}$ is a semisimple invertible skew-symmetric derivation of $\fg$ leaving $\fa^*$ invariant, then, according to Theorem \ref{T*down} and Remark 3, we have a semisimple $D\in\mbox{Der}(\fa)$ such that
$\overline{D}(x+f)=Dx-Hx-f\circ D$ for all $x\in\fa,\, f\in\fa^*.$
Since $D$ leaves $\mathfrak{z}(\fa)$ invariant and it is semisimple and invertible, we may choose the element $e\in\mathfrak{z}(\fa)$ such that $e\not\in\mathfrak{z}(\fh)$ of the proof of (i) to be an eigenvector of $D$. Now it suffices to prove that if $\fa$ is constructed as above then $\fa _1^*$ is invariant by $\overline{D}.$ 

Recall that 
$\fa _1^*=\mbox{${\Bbb C}$-span}\{x_1^*,x_2^*,\cdots ,x_n^*,e\},$
where $\{x_1^*,x_2^*,\cdots ,x_n^*\}$ is a basis of $\mathfrak h^*.$
Clearly, $\overline{D}(e)\in\fa_1^*$ because $\overline{D}(e)=\alpha e+He$ and
$0=B(\overline{D}(e),e)=B(He,e)$
implies that the projection of $He$ onto $e^*$ vanishes.
On the other hand, for every $j\leq n$ one has
$$B(\overline{D}(x_j),e)=-B(x_j,\overline{D}(e))=-B(x_j,\alpha e+He)=-B(x_j,\alpha e)=0,$$
which also shows that the projection of $\overline{D}(x_j)$ onto $e^*$ is zero.

\begin{coro}\label{reduction}
Let ${\fa}$ be a  Lie algebra such that and  $\theta\in Z^2({\fa},{\fa}^*)$ be a cyclic cocycle.

If the $T^*$-extension $\fg=T^*_\theta (\fa)$ is irreducible, then there exists a Lie algebra ${\mathfrak b}$ and a cyclic $\tilde{\theta}\in Z^2({\mathfrak b},{\mathfrak b}^*)$ such that ${\mathfrak z}({\mathfrak b})\subset [{\mathfrak b},{\mathfrak b}]_{{\mathfrak b}}$ and 
$ \fg=T^*_{\tilde{\theta}} ({\mathfrak b}).$

Moreover, if $\overline{D}$ is a semisimple skew-symmetric invertible derivation of $\fg$ leaving $\fa^*$ invariant then $\mathfrak b$ may be chosen such that $\mathfrak b^*$ is also stable by $\overline{D}.$

\end{coro}
\dem We will prove the result by induction on $l=\dim {\mathfrak z}(\fa)-\dim \Bigl([\fa,\fa]_{\fa}\cap {\mathfrak z}(\fa)\Bigl)$. 

If $l=0$ then obviously ${\mathfrak z}({\mathfrak a})\subset [{\mathfrak a},{\mathfrak a}]_{{\mathfrak a}}$ and we can choose $\mathfrak b=\fa.$ Suppose $l>0$ and that the result is valid for all $k<l$   and consider that 
$l=\dim {\mathfrak z}(\fa)-\dim \Bigl([\fa,\fa]_{\fa}\cap {\mathfrak z}(\fa)\Bigl)$. We clearly have that $\fa$ is as in Proposition \ref{red} and hence there exists a Lie algebra $\fa _1$ and a cyclic $\theta _1\in  Z^2(\fa_1,\fa^*_1)$ verifying
 $\fg=T^*_{\theta _1} (\fa_1)\, , \,\dim {\mathfrak z}(\fa _1)\le \dim {\mathfrak z}(\fa)\, ,$ and also $ 
\dim \Bigl([\fa _1,\fa_1]_{\fa_1}\cap {\mathfrak z}(\fa _1)\Bigl)=\dim \Bigl([\fa,\fa]_{\fa}\cap {\mathfrak z}(\fa)\Bigl)+1.
$
 We then have
 $$\dim {\mathfrak z}(\fa_1)-\dim \Bigl([\fa_1,\fa_1]_{\fa_1}\cap {\mathfrak z}(\fa_1)\Bigl)\le \dim {\mathfrak z}(\fa)-\dim \Bigl([\fa,\fa]_{\fa}\cap {\mathfrak z}(\fa)\Bigl)-1\le l-1$$
 and therefore, by the induction hypothesis, there exists a Lie algebra ${\mathfrak b}$ and a cyclic $\tilde{\theta}\in Z^2({\mathfrak b},{\mathfrak b}^*)$ such that ${\mathfrak z}({\mathfrak b})\subset [{\mathfrak b},{\mathfrak b}]_{{\mathfrak b}}$ and 
$ \fg=T^*_{\theta _1} (\fa_1)=T^*_{\tilde{\theta}} ({\mathfrak b}).$ 

The second  part  is immediate from statement (ii) in the proposition above.\qed

\bigskip

The next result gives the complete classification up to isometric isomorphism of quadratic Lie algebras with dimension less or equal that 8 admitting a symplectic structure.

First note that, up to dimension 4,  a nilpotent quadratic Lie algebra must be abelian. Hence, it suffices   to classify the cases of dimension 6 and 8 and therefore, according to  Theorem \ref{T*down},  certain $T^*$-extensions of those Lie algebras of dimension 3 and 4 which admit invertible derivations. There are only 2 nonisomorphic Lie algebras of dimension 3: the abelian Lie algebra and the Heisenberg algebra ${\cal L}_2$, which is the algebra spanned by three elements $x_1,x_2,x_3$ with the only non-trivial bracket  $[x_1,x_2]=x_3$. In dimension 4 one finds three nilpotent Lie algebras: the abelian 4-dimensional Lie algebra, the direct sum ${\cal L}_2\oplus \Bbb C$ and the filiform Lie alebra ${\cal L}_3$, which admits a basis $x_1,x_2,x_3, x_4$ where the only non-trivial brackets are $[x_1,x_2]=x_3$ and $[x_1,x_3]=x_4.$  According to Corollary \ref{reduction} a quadratic symplectic Lie algebra is either reducible or a $T^*$-extension of a Lie algebra whose center is contained in its derived ideal. Hence, for dimension less than or equal to 8, irreducible quadratic symplectic Lie algebras are given by the non-isomorphic $T^*$-extensions of the filiform Lie algebras ${\cal L}_2$ and ${\cal L}_3$ by means of a 2-cocycle $\theta$ compatible with some invertible derivation (in the sense that the map $\Theta$ of Theorem \ref{T*down} is a coboundary). The scalar cohomology and the deriviations of ${\cal L}_2$ and ${\cal L}_3$ are well-known (see, for example, \cite{borde} \cite{g-h}) and one easily verifies that for ${\cal L}_2$ such compatibility condition leads to $\theta =0$ whereas for ${\cal L}_3$, the associated 3-cocycle $\tilde{\theta}\in Z^3({\cal L}_3,\Bbb C)$ given by $\tilde{\theta}(x,y,z)=\theta(x,y)(z)$  must actually be a 3-coboundary  and thus, by \cite[Proposition 3.1]{borde}, $T^*_\theta({\cal L}_3)$ is isometrically isomorphic to the trivial $T^*$-extension $T^*_0({\cal L}_3)$. Finally, since the only quadratic symplectic Lie algebras up to dimension 4 are the abelian Lie algebras ${\Bbb C}^2$ and  ${\Bbb C}^4$, the reducible case follows at once and we get the following:

\begin{theor}
If $\fg$ is a complex quadratic  Lie algebra which admits a symplectic structure and $\dim (\fg)\le 8$, then $\fg$  is isometrically isomorphic to the trivial $T^*$-extension of one of the Lie algebras ${\Bbb C}$, ${\Bbb C}^2$, ${\Bbb C}^3$, ${\Bbb C}^4$, ${\cal L}_2$, ${\cal L}_2\oplus {\Bbb C}$ or ${\cal L}_3$.
\end{theor}

\section{Double extension of  quadratic symplectic Lie algebras} 

In her unpublished thesis, A. Aubert \cite{aubert} proved that every quadratic symplectic Lie algebra may be obtained from an abelian quadratic symplectic Lie algebra by a sequence of (generalised) symplectic double extensions by a line or a plane (see \cite{dard}, \cite{m-r91} for the formal definition of {\it symplectic double extensions}). For the sake of completeness and since they will be used in Section 4, we  include in this section one of Aubert's results.  For our purposes concerning Manin triples of Section 5 it is more accurate to describe quadratic symplectic Lie algebras in terms of quadratic double extension instead of symplectic double extension and, hence, the statements below look slightly different from those of Aubert's. We have considered, however, that such a difference does not justify the inclusion of the proofs in these paper. 

The main tool for all the constructions of quadratic Lie algebras below is the so-called {\it quadratic double extension}, described as follows:

\begin{defi}{\em \cite{m-r-travaux}, \cite{m-r85} $\,$
Let $({\fg},B)$ be a quadratic Lie algebra. Let ${\mathfrak b}$ another Lie algebra and $\psi:{\mathfrak b} \rightarrow
\mbox{\rm Der}_a({\fg},B)$ a representation of ${\mathfrak b}$ by means of skew-symmetric derivations of $({\fg},B)$. Then, the map  $\phi:{\fg} \times {\fg} \rightarrow
{{\mathfrak b}}^*$ defined by
$\phi(x,y)(z)=B(\psi(z)(x),y),$ for all $x,y\in {\fg}, \,  z\in \mathfrak b,$  
turns out to be a 2-cocycle for the trivial representation of 
${\fg}$ in ${\mathfrak b}^*$.
We now consider the central extension ${\mathfrak
b}^*\times_\phi {\fg}$ of $\fg$ by ${\mathfrak
b}^*$  by means of the 2-cocycle $\phi$ and the linear mapping ${\t}: {\mathfrak b}
\rightarrow {\mathfrak {gl}}({\mathfrak b}^*\times_\phi {\fg}) $ defined by:
$\t(x)_{{\vert} {{\mathfrak b}^*}}= \pi(x) \, , \,  \t(x)_{{\vert}{{\fg}}}\psi(x), $ for all $x\in {\mathfrak b},$
where  $\pi$ stands for the coadjoint representation of    ${\mathfrak b}.$  It follows that $\t$ is a representation of
${\mathfrak b}$ in ${\mathfrak b}^*\times_\phi {\fg} $ such that $\t(x) \in \mbox{\rm Der}({\mathfrak
b}^*\times_\phi {\fg}).$ 

 The Lie algebra $\tilde {\mathfrak g}$ given by the semi-direct product of  ${\mathfrak
b}^*\times_\phi {\fg}$  by
${\mathfrak b}$ via the representation $\t$ is called the {\it double extension of ${\fg}$
by ${\mathfrak b}$ by means of $\psi$}. 
}\end{defi}

If we identify the underlying vector space of $\tilde {\mathfrak g}$ with the direct sum of vector spaces  ${\mathfrak
b}\oplus {\fg} \oplus {\mathfrak b}^*$  then the bracket in $\tilde {\mathfrak g}$ is given by
\begin{eqnarray*}
 & & [y_1+x_1+f_1,y_2+x_2+f_2]= \Bigl([y_1,y_2]_{{\mathfrak b}}\Bigl)+\Bigl([x_1,x_2]_{{\fg}} +\psi(y_1)(x_2)-\psi(y_2)(x_1)\Bigl)+\\& &\hspace{0.5cm}
 \Bigl(\pi(y_1)(f_2)-\pi(y_2)(f_1)+\phi(x_1,x_2)\Bigl), \end{eqnarray*}
for $y_1,y_2\in {\mathfrak b},f_1,x_1,x_2\in {\fg},f_2\in {\mathfrak
b}^*.$ 
One can easily verify that the bilinear form $T$ given by
$T(y_1+x_1+f_1,y_2+x_2+f_2)=B(x_1,x_2)+f_1(y_2)+f_2(y_1), $ 
for $y_1,y_2\in {\mathfrak
b}, x_1,x_2\in {\fg} , f_1,f_2\in {\mathfrak b}^*$
is a scalar product on $\tilde {\mathfrak g}$ and hence $(\tilde {\mathfrak g},B)$ becomes a quadratic Lie algebra.

\medskip

In particular, if ${\mathfrak b}$ is a one-dimensional Lie algebra ${\mathfrak b}={\Bbb K}e$, we say that  $\f g$ is the {\it double extension of ${\fa}$  by means of the derivation $\psi(1)= \d.$}  Since this notion will be continously used from now on, we recall that  in this case one may consider 
 $\tilde {\mathfrak g}= {\Bbb K}e\oplus {\mathfrak g}\oplus {\Bbb
K}e^*$ endowed with the bracket
$[\alpha e + x + \alpha 'e ^*,\gamma e  + y + \gamma 'e ^*]_{\tilde {\mathfrak g}}[X,Y]_{\mathfrak g} +
\alpha \d(y) - \gamma \d(x) + B(\d(x),y)e^*,$  and the scalar product
$T(\alpha e +x + \alpha 'e^*,\gamma e +y + \gamma 'e^*)= B(x,y) + \alpha \gamma '+ \gamma
\a',$ for all $x,y \in {\mathfrak g}, \ \a, \a',\g, \gamma '\in {\Bbb K}$.

\medskip

\begin{lema}\label {LI1}  Let  $\mathfrak g$ be a Lie algebra, $B$  an invariant scalar
product on $\mathfrak g$ and  $D$ an invertible 
  skew-symmetric derivation of $\mathfrak g$.   Suppose that there exist $\delta \in Der_a({\mathfrak g},B)$,
$\lambda
\in {\Bbb K}\x \{0\}$ and $c\in {\mathfrak g}$ such that $[\d,D]-\lambda \delta = \adg (c)$ and let
$(\tilde{\fg} ,T)$ be the quadratic double extension of $({\fg} ,B)$ by means of $\d$.

 The linear endomorphism $\tilde D$ of  $\tilde {\mathfrak g}$ defined by
\begin{eqnarray*} & {\tilde  D}_{\vert_ {\mathfrak g}}= D + B(c,.) e^*, ~{\tilde  D}(e^*)= \lambda
e^*,~ 
 {\tilde  D}(e)= -\lambda e - c,& 
\end{eqnarray*}
is an invertible  derivation of $\tilde {\mathfrak g}$ which is
  skew-symmetric with respect to $T$ . 
\end{lema}
 
\begin{rema}  {\rm Obvously, if we set $\tilde{\o}(X,Y)=T(\tilde{D}X,Y)$, for $X,Y\in\tilde{\fg}$, then
$(\tilde {\mathfrak g},T,\tilde{\o}) $ is a quadratic symplectic Lie algebra which we will call
the {\it  quadratic symplectic double extension of
$(\mathfrak g ,B, \o)$ by the one-dimensional algebra by means of $(\d,c)$.} 
Actually, one can  easily see that  $(\tilde {\mathfrak g},\tilde{\o}) $ is the symplectic
double extension of 
  $(\mathfrak g , \o)$ by means of  $\d$ and $z= D^{-1}(\d(c))$
(see \cite{dard}  ou \cite {m-r91}  for the notion of symplectic double extension). }\end{rema}

 \begin{theor}\label {TI1} Let $(\tilde {\mathfrak g},T,\tilde{\o}) $ be a quadratic symplectic Lie
algebra over a field $\Bbb K$ and let
$\tilde{D}$ be the invertible derivation such that $\tilde{\o}(X,Y)=T(\tilde{D}X,Y)$ for every
$X,Y\in\tilde {\mathfrak g}$.
\begin{enumerate} 
\item[(i)] If there exists an one-dimensional central ideal of $\tilde {\mathfrak g}$ which
remains invariant by the derivation $\tilde D$, then  $(\tilde {\mathfrak g},T,\tilde{\o}) $ is 
the quadratic symplectic double extension of a quadratic symplectic Lie algebra
$(\mathfrak g ,B, \o)$ by the one-dimensional algebra by means of a pair $(\d,c)$.
\item[(ii)] In particular, if  $\Bbb K$ is algebraically closed, then $(\tilde {\mathfrak g},T,\tilde{\o}) $
may be obtained from the two-dimensional abelian Lie algebra by a
sequence of quadratic symplectic double extensions by the one-dimensional algebra
 \end{enumerate}
\end{theor}
\section{Manin algebras  and quadratic symplectic Lie algebras}

A particular and very interesting case of quadratic double extension concerns Manin triples. Using this procedure developped in \cite{m-r88} and the results of section 4 we will describe more precisely the structure of quadratic symplectic Lie algebras.

\begin{defi} {\rm Let $\f g$ be   Lie algebra and let $\f U$, $\f V$ be two Lie 
subalgebras of $\f g$. The triple $({\f g},{\f U},{\f V})$   is called a {\it Manin
triple} if  $\fg$ is the direct sum of vector subspaces ${\f g}= {\f U}\oplus{\f V}~$ and there exists    an invariant scalar product $B$ on
$\f g$  such that ${\f U}$ and ${\f V}$ are  completely isotropic with respect to $B$. In such case we will also say that
$({\f g}= {\f U}\oplus{\f V},B) $ is a {\it Manin algebra } and  ${\f g}= {\f U}\oplus{\f V}$  
a {\it Manin decomposition} of $\f g.$
}
\end{defi}

\begin{rema} {\em
\begin{enumerate}\item[$\,$]
\item Obviously, any trivial $T^*$-extension $\fg=T^*_0\fa$ provides a Manin triple $(\fg, \fa^*,\fa).$
\item It is well known  that there is one-to-one  correspondence between Lie bialgebra structures on $\f
U$  and Manin triples $({\f g},{\f U},{\f V}) $ \cite{CP}. 
\end{enumerate}}
\end{rema}

\begin{theor}\label{alberto} {\em \cite{m-r88}} Let $({\f G}= {\f U}\oplus{\f V},B)$ be  a Manin algebra. 
 Let $\delta$ be a skew-symmetric derivation of $(\f G,B)$ such that $\delta({\f V})\subseteq {\f V}.$ 

The Lie algebra obtained by quadratic double extension of $\f
G$  by means of $\delta$ is also a  Manin algebra $({\f g}= {\f U}'\oplus{\f V}',{\tilde B}), $
where ${\f U}'$ is the direct sum of Lie algebras
${\f U}'= {\f U}\oplus {\Bbb K}e^*$ and  ${\f V}'= {\f V}\oplus {\Bbb K}e$ with the Lie structure given by the semidirect product by means of $\delta$. \end{theor}

\begin{defi} {\rm The Manin algebra $ ({\f g}= {\f U}'\oplus {\f V}',{\tilde B}) $ obtained in the theorem above is called the
{\it double extension of  the Manin algebra $ ({\f G}= {\f U}\oplus {\f V},B)$ by the one-dimensional
Lie algebra (by means of $\delta$).}}
\end{defi}

\begin{prop}\label {PI1} Let  $({\f  g}= {\f U'}\oplus{\f V'},B')$ be  a  Manin algebra of dimension
$n$. 

If  either ${\f z}({\f g})\cap {\f U'}\not= \{0\} ~\mbox{or}~ {\f z}({\f g})\cap {\f
V'}\not= \{0\},$ then $({\f  g}= {\f U'}\oplus{\f V'},B')$ is a double extension of a Manin algebra
$({\f  G}= {\f U}\oplus{\f V},B)$ of dimension $n-2$ by the one-dimensional Lie algebra.
\end{prop}

\dem  If ${\f z}({\f g})\cap {\f U'}\not= \{0\},$ then there exists $e^* \in {\f z}({\f
g})\x\{0\}$ such that ${\f I}={\Bbb K}e^* \subseteq {\f U'}$   and, hence,  there exists $e \in {\f
V'}$ such that $B'(e^*,e)= 1.$ Consequently, if ${\f G}= ({\Bbb K}e^*\oplus {\Bbb K}e)^{\perp}$ denotes the orthogonal of ${\Bbb K}e^*\oplus
{\Bbb K}e$ with respect to $B'$, one easily verifies that $B:= B'_{\vert_{{\f G}\times {\f G}}}$ is an invariant scalar
product on $\f G$ and that $(\fg, B')$ is the quadratic double extension of $(\f G, B)$ by means of the derivation $\delta=\mbox{ad}_\fg (e)_{\vert\f G}.$ 

Now, 
 ${\f U}= {\f G}\cap {\f U'}$ is a Lie subalgebra of $\f G$ and 
${\f U'}= {\Bbb K}e^* \oplus {\f U}$  because $\f U'$ is a subalgebra of $\f g$ and $e^* \in {\f U'}.$
Moreover, $B({\f U},{\f U})= B'({\f U},{\f U})=\{0\}.$ The fact that ${\f g}= {\f U'} \oplus {\f
V'} ~\mbox{and}~ {\f U}'\subseteq {\f I}^{\perp}$ implies that ${\f I}^{\perp}= {\f U'} \oplus ({\f
V'}\cap {\f I}^{\perp}).$ Since $B'({\f V'}\cap {\f I}^{\perp},{\Bbb K}e^*\oplus {\Bbb K}e)= \{0\},$
we immediately get that 
${\f V'}\cap {\f I}^{\perp}\subseteq {\f G}$ and, furhter, ${\f G}= ({\f U'}\cap {\f G}) \oplus
({\f V'}\cap {\f I}^{\perp})$ because ${\f G}\subseteq {\f I}^{\perp}.$ Consequently, ${\f G}= {\f
U}\oplus {\f V}$ where ${\f V}= {\f V'}\cap {\f I}^{\perp}$ is a Lie subalgebra of $\f G$ such
that
$B({\f V},{\f V})= \{0\}.$ This  proves that $({\f G}= {\f U}\oplus {\f V}, B)$ is a
Manin algebra of $\mbox{dim}{(\f G)}= \mbox{dim}{(\f g)}-2.$ Finally, since $e \in {\f V'},$ then $\delta ({\f V})=\mbox{ad}_{\f
g}(e)({\f V})\subseteq {\f V}$ and we conclude that ${\f g}= {\f U'}\oplus {\f
V'}$ is the double extension of the Manin algebra  ${\f G}= {\f U}\oplus {\f
V}$ by  means of $\delta$.\qed           

\begin{theor}\label{nilman} Let  $({\f  g}= {\f U'}\oplus{\f V'},B')$ be  a non-zero Manin algebra. If
$\f g$ is a  nilpotent Lie algebra, then ${\f z}({\f g})\cap {\f U'}\not\{0\}~\mbox{and}~{\f z}({\f g})\cap {\f V'}\not= \{0\}.$ 
As a consequence, $({\f  g}= {\f U'}\oplus{\f V'},B')$ is a double extension of a Manin algebra of dimension $\mbox{\rm dim} (\fg)-2$.   
\end{theor}

\dem   Since
$\f g~$ and $~\f U'$ are non-zero  nilpotent Lie algebras then ${\f z}({\f g})\not= \{0\}
$ {and} ${\f z}({\f U'})\not= \{0\}.$   Suppose that ${\f z}({\f g})\cap {\f U'}= \{0\}.$    We
set ${\f L}_0= [{\f V'},{\f z}({\f U'})]$, it is clear that ${\f L}_0\not=\{0\}~\mbox{and}~ [{\f
U'},{\f L}_0]\subseteq{\f L}_0.$  Let us consider
$v\in {\f V'}$ and $u\in {\f z}({\f U'}) $ and write $[v,u]= x + y$ where $x\in {\f U'}$ and $y\in {\f V'}.$
Now, if $t\in {\f U'},$ then $B(y,t)= B'([v,u],t)= B'(v,[u,t])= 0$ and the non-denereracy of $B'$ shows that $y=0.$ It
follows that $ {\f L}_0\subseteq {\f U'}.$ Therefore, ${\f L}_0$ is a non-zero ideal of
$\f U',$ which implies that ${\f L}_0\cap {\f z}({\f U'})\not= \{0\}$ because $\f U'$ is a nilpotent
Lie algebra.

Now, consider ${\f L}_1= [{\f V'},{\f L}_0\cap {\f z}({\f U'})].$ The fact that ${\f L}_0\cap {\f
z}({\f U'})\not= \{0\}$ and ${\f z}({\f g})\cap {\f U'}= \{0\}$ implies that ${\f L}_1\not=\{0\}.$
Further, as $B'$ is an invariant scalar product on $\f g$, we have  that ${\f L}_1\subseteq {\f U'}.$
Since
$[{\f U'},{\f L}_1]= [{\f L}_0\cap {\f z}({\f U'}),[{\f U'},{\f V'}]]\subseteq [{\f L}_0\cap {\f z}({\f
U'}),{\f V'}]= {\f L}_1,$ then ${\f L}_1$ is an ideal of $\f U'.$ Consequently, ${\f L}_1\cap {\f
z}({\f U'})\not=\{0\}.$ If we  repeat  this process sucessively, we get a sequence $({\f L}_n)_{n \in {\Bbb N}}$ of non-zero
ideals of $\f U'$ defined by ${\f L}_0= [{\f V'},{\f z}({\f U'})] ~\mbox{and} ~ {\f L}_n= [{\f V'},{\f
L}_{n-1}\cap {\f z}({\f U'})],$ for $n\geq 1.$ Further, it is obviously verified that  ${\f
L}_n\subseteq {\cal C}^{n+1}({\f g}),~$ where $({\cal C}^{n}({\f g}))_{n \in {\Bbb N}}$ stands for the central descending series 
of $\f g.$ The fact that $\f g\not= \{0\}$ is a nilpotent Lie algebra implies that there exists $k \in {\Bbb N}$ such
that ${\cal C}^{k}({\f g})=\{0\}$. Consequently, $ {\f L}_{k-1}=\{0\}$
which contradicts the fact that ${\f L}_n\not=\{0\}$ for all $n \in {\Bbb N}.$ We conclude that ${\f
z}({\f g})\cap {\f U'}\not= \{0\}.$ The same reasoning shows that ${\f z}({\f g})\cap {\f
V'}\not= \{0\}$.
Now, the final part of the statement is immediate from Proposition \ref{PI1}.\qed

\begin{coro} Every nilpotent Manin algebra can be obtained by a sequence of quadratic double extensions by  one-dimensional Lie algebras starting from a 2-dimensional abelian Manin algebra. 
\end{coro}

\begin{defi} {\rm We will say that a  Manin algebra $({\f  g}= {\f U}\oplus{\f V},B) $  is a {\it special symplectic Manin algebra} if there exists   a symplectic structure $\omega$ on the Lie algebra $\f g$
such that 
$\omega({\f U},{\f U})= \omega({\f V},{\f V})= \{0\}.$}
\end{defi}

The following lemma follows easily:

\begin{lema} \label {SIA} A Manin  algebra $({\f  g}= {\f U}\oplus{\f V},B)$ is  special symplectic if and
only if there exists an invertible skew-symmetric derivation $D$ of $(\f g, B)$ such that
$D({\f U})\subseteq {\f U}~\mbox{and}~D({\f V})\subseteq {\f V}.$ 
\end{lema}

\noindent {\bf Notation } If a Manin algebra $({\f  g}= {\f U}\oplus{\f V},B) $ admits a special symplectic structure we will note the corresponding symplectic quadratic algebra by either $({\f  g}= {\f
U}\oplus{\f V},B,\omega)$ {or} by  $({\f  g}= {\f U}\oplus{\f V},B,D),$ where $D$ denotes the derivation given in the lemma.  

\medskip

Our main interest in the study of special symplectic Manin algebras lies in the fact that every quadratic symplectic Lie algebra is so. More precisely:

\begin{prop} Let $({\f g},B,\omega)$ be a quadratic symplectic Lie algebra over an algebraically closed field $\Bbb K$. There exist
two Lie subalgebras   $\f U$ and $\f V$    of $\f g$ such that $({\f  g}= {\f U}\oplus{\f V},B,\omega)$ 
is a special symplectic Manin algebra.
\end{prop}
\dem  Consider  $D \in \mbox{Der}_a({\f g})$ invertible such that $\omega(x,y)= B(D(x),y)$  for all $x,y
\in {\f g}.$ The proof follows the same argument used in \cite[Corollary 2.16]{al-mi}. Let
$\mbox{\rm Sp}(D)$ denote the set of all distinct eigenvalues of
$D$. Since every field of characteristic zero is a $\Bbb Q$-vector space, we can consider the $\Bbb Q$-linear span $\mathcal S$ of $\mbox{Sp} (D)$ (notice that $\mbox{dim}_{\Bbb Q}\ {\mathcal S}\leq 
\mbox{dim}_{\Bbb K}\ {\mathcal \fg}/2$ since $D$ is skew-symmetric). Let us fix a basis ${\cal H}=\{\alpha_1,\alpha_2,\cdots,\alpha_k\}$ of $\mathcal S$ and let us consider on $\mathcal S$ the lexicographic order with respect to ${\cal H}$, this is to say:
$\sum_{i=1}^km_i\alpha_i >0$ if and only if there exists $i_0\leq k$ such that $m_{i_0}>0 \mbox{ and } m_j=0 \mbox{ for }j<i_0.$
Obviouosly, with this order we have that 
$\alpha +\beta >0  $ if $\alpha,\beta>0$.  Now, define $\mbox{\rm Sp}^+$ and $\mbox{\rm Sp}^-$ respectively as the sets of positive and negative elements of $\mbox{\rm Sp}(D)$.
It is clear that  $\mbox{\rm Sp}(D)={\rm Sp}^+\cup {\rm
Sp}^- \ \mbox{and } {\rm Sp}^+\cap {\rm Sp}^-= \emptyset.$ For each $\lambda \in {\rm
Sp}(D),$ let us consider  $\fg (\lambda )=\{ x\in\fg
\, : \,  (D-\lambda \mbox{\rm id}_\fg)^{\mbox{\rm\small dim}(\fg)}(x)=0\}. $   The fact that $D$ is skew-symmetric with respect to $B$ implies that $B({\f g}(\lambda
),{\f g}(\lambda'))= \{0\}$ for all $\lambda,\lambda'\in {\rm Sp}(D)$ such that
$\lambda+\lambda'\not=0.$ Now, if we set $${\f U}=\sum_{\lambda\in\mbox{\rm Sp}^+} \fg
(\lambda )\, ,\quad {\f V}=\sum_{\lambda\in\mbox{\rm Sp}^-} \fg (\lambda ),$$
then $\f U$ and $\f V$ are two completely isotropic Lie subalgebras of $\f g$ which are stable by $D$.   We conclude that $({\f  g}= {\f U}\oplus{\f V},B,\omega)$ is a
special symplectic Manin algebra.\qed  
 
\medskip

The inductive description of special symplectic Manin algebras  may also be done by a double extension procedure using only Manin algebras. In order to show this, we begin with the following analog of  Lemma \ref{LI1}:  

\begin{lema} Let $({\f  G}= {\f U}\oplus{\f V},B,D)$ be a special symplectic Manin algebra and suppose that there exist $\delta
\in \mbox{\rm Der}_a({\mathfrak G},B)$, $\lambda \in {\Bbb K}\x \{0\}$ and $c\in {\mathfrak V}$ such that
$[\d,D]-\lambda \delta = \ad_{\f G} (c)$ and $\delta({\f V})\subseteq {\f V}.$

 Let $({\f  g}= {\f U}'\oplus{\f V}',{B'})$  be the double
extension of the Manin algebra    $({\f  G}= {\f U}\oplus{\f V}, B)$ by means the derivation $\delta$ (constructed as in Theorem
\ref{alberto}).

The skew-symmetric invertible derivation ${\tilde D}$ of $(\f g,B')$ defined by:
\begin{eqnarray*} & {\tilde  D}_{\vert_ {\mathfrak G}}= D + B(c,.) e^*, ~{\tilde  D}(e^*)= \lambda
e^*,~ 
 {\tilde  D}(e)= -\lambda e - c,& 
\end{eqnarray*}
verifies ${\tilde  D}({\f U}')\subseteq {\f U}'~\mbox{and}~ {\tilde  D}({\f V}')\subseteq {\f
V}'.$ Consequently, $({\f  g}= {\f U}'\oplus{\f V}',{ B'},{\tilde  D})$ is a also a special symplectic
Manin algebra.\end{lema}
\dem Bearing in mind Lemma \ref{LI1}, it suffices to prove that $\f U$ and $\f V$ are stable by ${\tilde  D}.$ Since $D({\f U})\subseteq {\f U}$ and ${\tilde  D}(e^*)\lambda e^*,$ we immediately get ${\tilde  D}({\f U}')\subseteq {\f U}'.$ The fact that $c\in {\f V}$ and
$B({\f V},{\f V})= \{0\}$ implies that ${\tilde  D}({\f V})= D({\f V})\subseteq {\f V}.$
Consequently ${\tilde  D}({\f V}')\subseteq {\f V}'$ because ${\tilde  D}(e)= -\lambda e - c.$ \qed

\bigskip 

The special symplectic Manin algebra $({\f  g}= {\f U}'\oplus{\f V}',B',{\tilde D})$ of the Lemma will be called the {\it special double extension of the symplectic Manin algebra $({\f 
G}= {\f U}\oplus{\f V},B,D)$ by means of $(\delta,c).$}

\smallskip

\begin{theor}\label {Classif}
Let $({\f  g}= {\f U'}\oplus{\f V'},B',\tilde{D})$ be a special symplectic Manin algebra over  $\Bbb K$.
\begin{enumerate} 
\item[(i)] If $\tilde D$ admits an eigenvector $z\in {\f z}({\f g})\cap {\f U'}+{\f z}({\f g})\cap {\f V'}$, then  $({\f  g}= {\f U'}\oplus{\f V'},B',\tilde{D})$ is a  special double extension of a special
symplectic Manin algebra $({\f  G}= {\f U}\oplus{\f V},B,D)$ by means of a pair $(\delta,c)$ where
$c\in {\f V}~\mbox{and}~ \delta$ is a skew-symmetric derivation of $\f G$ leaving  ${\f V}$ invariant.
\item[(ii)] In particular, if  $\Bbb K$ is algebraically closed, then $({\f  g}= {\f U'}\oplus{\f V'},B',\tilde{D})$ may be obtained from a two-dimensional special symplectic Manin algebra  by a
sequence of special double extensions by the one-dimensional algebra.
\end{enumerate}\end{theor}
\dem 
Let us consider $z\in {\f z}({\f g})\cap {\f U'}+{\f z}({\f g})\cap {\f V'}$, $z\neq 0$ such that $\tilde{D}z=\lambda z$ for some $\l \in {\Bbb K}\x \{0\}$ and put $z=u+v$ where $u\in {\f z}({\f g})\cap {\f U'}, \ v\in{\f z}({\f g})\cap {\f V'}.$ Since both $U'$ and $V'$ are stable by $\tilde{D}$, it is clear that $\tilde{D}u=\lambda u$ and $\tilde{D}v=\lambda v$ and either $u\neq 0$ or $v\neq 0$. We will suppose that $u\neq 0$ (otherwise, we can change the roles of $U'$ and $V'$). If we set $e^*=u$ then we have, as in Proposition   \ref {PI1} that $({\f  g}= {\f U'}\oplus{\f V'},B')$ is the double extension of a Manin algebra $({\f  G}= {\f U}\oplus{\f V},B)$ where $\f U'={\f U}\oplus {\Bbb K}e^*$, $\f V'={\f V}\oplus {\Bbb K}e$ and $B'(e^*,e)=1$ by means of the derivation $\delta =\mbox{ad}_\fg (e)_{\vert_{\f G}}.$

Since the ideal $\f I= {\Bbb K}e^*$ is invariant by the skew-symmetric mapping $\tilde{D}$, so is its orthogonal ${\f I}^\bot={\Bbb K}e^*\oplus\f G$. Now, if $p: {\f I} ^\bot={\Bbb K}e^*\oplus\f G\rightarrow \f G$ denotes  the projection $p(\alpha e^*+x)=x$ for $\alpha \in {\Bbb K}, x\in\f G$, then one can easily verify that $D=p\circ \tilde{D}_{\vert\f G}  $ is an invertible skew-symmetric derivation of $(\f G, B)$. Further, as ${\f U}={\f U'}\cap \f G$ and ${\f U'}$ is invariant by $\tilde{D}$, we have that $D(\f U)\subset \f U$ and, analogously, we prove $D(\f V)\subset \f V$. Since $\tilde{D}$ is skew-symmetric, $B(e^*,e)=1$ and  $\tilde{D}(e)\in\f V'$, one immediately obtains that there exists $c\in V$ such that  $ {\tilde  D}(e)= -\lambda e - c$ and ${\tilde  D}_{\vert_ {\mathfrak G}}= D + B(c,.) e^*.$ 
Further, from the fact that $\tilde{D}$ is a derivation, one easily deduces that 
$ [\d,D']-\lambda \delta = \ad_{\f G}(c)$ and, therefore, $({\f g}= {\f U'}\oplus{\f V'},B',\tilde{D})$ is the  special double extension of $({\f  G}= {\f U}\oplus{\f V},B,D)$ by means of $(\delta,c)$ as claimed in (i).

Now, if $\Bbb K$ is algebraically closed, we can always find an eigenvector of $\tilde D$ in ${\f z}({\f g})\cap {\f U'}$ since, according to Theorem \ref{nilman}, it is a non-zero subspace and it is obviously stable by all derivations of $\fg$ which leave $ \f U'$ invariant. If we apply (i) sucessively then (ii) follows.\qed

%%%%%%%%%%%%%%%%%%%%%%%%%%%%
%%%%%%%% REFERENCES
%%%%%%%%%%%%%%%%%%%%%%%%%%%%

\vskip 0,5 cm

\n {\it Authors' addresses:}
\\
\n {\bf I. Bajo}: Depto. Matem\'atica Aplicada II, E.T.S.I.
Telecomunicaci\'on, Universidad de Vigo, 36280 Vigo (Spain); ibajo@dma.uvigo.es \\
\n {\bf S. Benayadi}: Universit\'e Paul Verlaine-Metz, LMAM CNRS UMR 7122, Ile du Saulcy, 57045 Metz cedex (France); benayadi@poncelet.univ-metz.fr\\
\n {\bf A. Medina}: Dept. Math\'ematiques, CC 051, Universit\'e de Montepellier II, Place Eug\`ene Bataillon, 34095 Montpellier cedex 5 (France); medina@math.univ-montp2.fr
\end{document}